\documentclass[11pt,english]{amsart}
\def\margin_comment#1{\marginpar{\sffamily{\tiny #1\par}\normalfont}}

\usepackage{amssymb,euscript}

\usepackage{amsfonts}
\usepackage{geometry}
\usepackage{babel}

\usepackage{setspace}
\usepackage{color}

\usepackage{amsmath,amssymb,amsthm}
\usepackage{epsfig}

\usepackage{tikz}
 \usepackage{amsmath}
\usepackage{tikz-cd}
\usetikzlibrary{matrix, calc, arrows}
\usepackage{tikz}
\usetikzlibrary{arrows,chains,matrix,positioning,scopes,}
\usepackage{lipsum}
\usepackage{mathtools}
\usepackage{mathrsfs}
\usepackage{graphicx}
\usepackage{graphics}
\usepackage{setspace}
	\usepackage{color}

\usepackage{tikz}
\usetikzlibrary{trees}
\usepackage{tikz-cd}
\usetikzlibrary{fit,shapes,automata,backgrounds,petri,positioning,scopes,matrix, calc, plothandlers,plotmarks,patterns}

\usepackage{enumerate}
\usepackage{float}

\usepackage{tikz-3dplot}
\usepackage{pgfplots}
\usepackage{tkz-graph}
\usepackage{tikzscale}
\GraphInit[vstyle = Normal]
\tikzset{
  LabelStyle/.style = {minimum width = 2em, 
                        text = red, font = \bfseries },
  VertexStyle/.append style = { inner sep=2pt,
                                font = \Large\bfseries, fill},
  EdgeStyle/.append style = {->, bend left} }

\makeatletter
\tikzset{join/.code=\tikzset{after node path={%
\ifx\tikzchainprevious\pgfutil@empty\else(\tikzchainprevious)%
edge[every join]#1(\tikzchaincurrent)\fi}}}
\makeatother
\tikzset{>=stealth',every on chain/.append style={join},
         every join/.style={->}}
\tikzstyle{labeled}=[execute at begin node=$\scriptstyle,
   execute at end node=$]
\newtheorem{thm}{Theorem}[section]
\numberwithin{equation}{section} 
\numberwithin{figure}{section} 
\theoremstyle{plain}
\newtheorem*{thm*}{Theorem}
\theoremstyle{definition}
\theoremstyle{plain}
\newtheorem{thm_A}{Theorem}
\newtheorem*{defn*}{Definition}

\theoremstyle{plain}

\theoremstyle{plain} 

\theoremstyle{plain}

\theoremstyle{remark}
\newtheorem{ex}[thm]{Example}
\theoremstyle{remark}

\theoremstyle{plain}

\theoremstyle{plain}

\theoremstyle{plain}
\newtheorem{lem}[thm]{Lemma} 
\theoremstyle{definition}
\newtheorem{defn}[thm]{Definition}
\newtheorem*{acknowledgment}{Acknowledgment}
\newtheorem*{acknowledgment*}{Addentum}

\theoremstyle{plain}
\newtheorem*{ex*}{Example}
\theoremstyle{plain}
\newcommand{\Div}[1]{\operatorname{Div}(#1)}

\AtBeginDocument{
  
}

\begin{document}

\title[A note on Garside monoids and  $\mathscr{M}$-Braces ]{A note on Garside monoids  and  $\mathscr{M}$-Braces}
\author{Fabienne Chouraqui}
\begin{abstract}
	We define an algebraic structure  similar to that of a semiring, but  without some of the requirements. As it is somehow also similar to the structure of  left brace, we call it  an  $\mathscr{M}$-brace.
	We present a connection between Garside monoids and more generally lcm-monoids with  this  algebraic structure.
An  lcm-monoid  $M$  is a left-cancellative  monoid  such that $1$ is the unique invertible element in $M$,   and every pair of elements in $M$ admit an lcm  with respect to left-divisibility.  The class of lcm-monoids contains the  Gaussian, quasi-Garside and Garside monoids. 
\end{abstract}
\maketitle
\section*{Introduction}
	A left brace is a triple $(\mathcal{B},+,\cdot)$, where  $(\mathcal{B},+)$ is an abelian group, 
$(\mathcal{B},\cdot)$ is a group, and there is a left-distributivity-like  axiom that relates between the two operations in $\mathcal{B}$ \cite{rump_braces}.   In analogy with a left brace, we define a left  $\mathscr{M}$-brace to be a triple $(\mathcal{B},+,\cdot)$, where  $(\mathcal{B},+)$ is a commutative monoid, 
$(\mathcal{B},\cdot)$ is a monoid, and the  axiom of left  distributivity holds.  This  algebraic structure  is similar to that of a semiring, but  without the requirement of  the following two  axioms of a  semiring:  the right-distributivity axiom and  the product of  the zero element in $(\mathcal{B},+)$ with any element in $\mathcal{B}$ is the zero element.  In this short note, we show there is a  connection between lcm-monoids and left  $\mathscr{M}$-braces. The approach  is constructive and  since   Garside monoids, quasi-Garside monoids  and  Gaussian monoids  are subclasses of  lcm-monoids, this gives a connection between  these intermediate, better known, classes and  left  $\mathscr{M}$-braces.  More formally, we show:
	\begin{thm_A}\label{thm1}
		Let $M $ be an  lcm-monoid. Then there exists a left   $\mathscr{M}$-brace  $(\mathcal{B},\oplus,\cdot)$ such that $(\mathcal{B},\cdot)$ is isomorphic to $M$.
		Furthermore, if  $M $  is  a Gaussian monoid, then there exists a left   $\mathscr{M}$-brace  $(\mathcal{B},\oplus,\cdot)$,  and a right $\mathscr{M}$-brace $(\mathcal{B},+,\cdot)$,  such that $(\mathcal{B},\cdot)$ is isomorphic to $M$.
	\end{thm_A}
 In \cite{chou-F},   we define a left partial brace to be a triple $(\mathcal{B},\oplus,\cdot)$ such that $(\mathcal{B},\oplus)$ is a commutative partial  monoid (as defined in \cite{partial-s}),  $(\mathcal{B},\cdot)$ is an inverse monoid,  and the 	left-distributivity axiom holds. This definition of partial left brace enables us to consider Gaussian groups with a complemented presentation which satisfies a  certain condition and to  show there is a connection between them and left partial braces. 
 	\margin_comment{\textcolor{red}{added defn}} A \emph{complemented presentation} of a monoid $M$  is  $\operatorname{Mon}\langle X \mid R_{\lambda}\rangle$, where $X$ is a set, $X^*$ is  the free monoid generated by $X$ and  $R_{\lambda}$ is the set of defining relations $\{x\lambda(x,y)=y\lambda(y,x)\mid x,y \in X\}$, where  $\lambda: X\times X \rightarrow X^*$ is a  function satisfying $\lambda(x,x)=1$ for every $x \in X$ \cite{dehornoy}.

		\begin{thm_A}\label{thm2}
			Let $X$ be a set and $X^*$ be the free monoid generated by $X$. Let $M$ be a Gaussian monoid  with a complemented presentation  $M=\operatorname{Mon}\langle X \mid R_{\lambda}\rangle$, such that for every $x \in X$, the functions $\lambda(x,-):X \rightarrow X^*$ are injective. Let  $G$ be the  group of fractions of $M$.	Then there exists a partial  left   brace  $(\mathcal{B},\oplus,\cdot)$  such that $(\mathcal{B},\cdot)$ is isomorphic to $G$.
	\end{thm_A}

\margin_comment{\textcolor{red}{Theorem 2:
		added condition}}
The paper is organized as follows. In Section $1$, we give some preliminaries on braces, on monoids,  on lcm-monoids and
  Garside monoids.  In Section $2$, we prove the main results.
\section{Preliminaries}

\subsection{Preliminaries on Braces}
In \cite{rump_braces}, W. Rump introduced braces as a generalization of radical rings related with solutions of the Yang-Baxter equation.   In subsequent papers, he developed the theory of this new algebraic structure. In \cite{brace}, the authors give another equivalent definition of a brace and study its structure. 
\begin{defn} \cite{brace}
	A \emph{left brace}  is a set  $G$ with two operations, $+$ and $\cdot$, such that $(G,+)$ is an abelian group, $(G,\cdot)$ is a group and for every $a,b,c \in G$:
	\begin{equation}\label{eqn-left-dist-like}
a \cdot (b+c) = a \cdot b+a \cdot c -a
	\end{equation}
	The groups  $(G,+)$  and $(G,\cdot)$ are called \emph{the additive group} and \emph{the multiplicative group of the brace},  respectively.\\
	A right brace is defined similarly,  replacing condition  \ref{eqn-left-dist-like} by the following condition:
		\begin{equation}\label{eqn-right-dist-like}
(a+b) \cdot c  = a \cdot c+b \cdot c -c
	\end{equation}

\end{defn}
A \emph{two-sided brace} is a left and right brace, that is both $(*)$  and   $(**) $ are satisfied. From the definition of a left brace $G$, it  follows  that the  identity of the multiplicative group of $G$ is equal to the identity  of the additive group of $G$. Additionally,  for every $a,b,c \in G$, $a\cdot (b-c)=a\cdot b-a \cdot c +a$. 
\subsection{Preliminaries on  monoids and groups of fractions}
\label{sec_bcgd_monoids}

Let $X$ be a non-empty set. Let  $X^*$ denote   the free monoid generated  by $X$.  A \emph{monoid presentation}  is a pair $(X,R)$, where $R$ is a list of   pairs of words of  $X^*$. We denote by $\operatorname{Mon}\langle X\mid R\rangle$  the monoid $M$ presented by $(X,R)$,  that is, the monoid $ (X^*\big/\equiv R)$,  where $\equiv R$ is the congruence on $X^*$  generated by $R$.   A monoid $M$ is  \emph{left cancellative} if for every $a,b,b'$ in $M$, the equality~$ab=ab'$ in $M$ implies $b=b'$ in $M$, $M$ is  \emph{right cancellative} if for every $b,b',c $ in $M$, the equality~$bc=b'c$ in $M$ implies $b=b'$ in $M$, and $M$ is  \emph{cancellative}  if it is both left and right cancellative.

 A \emph{group of right (or left) fractions} of a monoid $M$ is a group $G$ such
 that $M$  embeds in $G$  and moreover for every $ g\in G$, there are $a,b \in M$ such that     $g=ab^{-1}$ (or there are $c,d \in M$ such that   $g=c^{-1}d$ ).  A commutative monoid embeds into an abelian group of fractions if and only if it is cancellative \cite{clifford}. 
 More generally, Ore's conditions state that a  cancellative monoid $M$  that satisfies  for every $a,b \in M$, $aM\cap bM \neq \emptyset$, (that is 
 for every $a,b \in M$, there are  $a',b' \in M$, such that $aa'=bb'$),  embeds into a group of
left  fractions.  And a  cancellative monoid $M$  satisfying for every $a,b \in M$, $Ma\cap Mb \neq \emptyset$, (that is 
 for every $a,b \in M$, there are $a'',b'' \in M$, such that $a''a=b''b$), embeds into a group of
right fractions.  Moreover, this  group of  left (or right)  fractions  is unique up to isomorphism \cite{ore }, \cite[p.300]{clifford}.   If both  Ore's conditions are satisfied, the groups of left and  right fractions  obtained coincide and we simply say the group of fractions:  for every $g \in G$, there are  $a,b,c,d \in M$, such that $g=ab^{-1}=c^{-1}d$.
 
  \begin{defn}
 	Let $M$ be a monoid.  Let $a,b,c,c',d$ be elements in $M$. 
 	\begin{enumerate}[(i)]
 		\item   $a$ is a \emph{left divisor  of  $d$ or left divides $d$}, if there is an element $a'$ in $M$ such that $d=aa'$.  
 		\item  $c$
 		is
 		\emph{a common multiple  with respect to left-divisibility} of  $a$ and $b$, if $a$ and $b$ are left divisors of $c$, that is $cM \subseteq  aM\cap bM$. 
 		\item  
$c$ is  \emph{a  least common multiple  with respect to left-divisibility } of $a$ and $b$, denoted by $c=a\vee b$,  if $c$  	is a common multiple of  $a$ and $b$,  and $c$ left divides any other  common multiple of  $a$ and $b$, that is $cM=aM\cap bM$. 
 	\end{enumerate}
 	\end{defn}
\begin{defn}
Let $M$ be a monoid and let $a,b$ be elements in $M$. Say that $d$
is \emph{a left greatest common divisor  of $a$ and $b$ with respect to left-divisibility} of $a$ and $b$,   if  $d$ is a  left divisor of $a$ and $b$, and additionally   if there is  $d' \in M$ such that $d'$ is a left divisor of $a$ and $b$, then
	$d'$ is  a left divisor of $d$.
\end{defn}
As it is common, we write \emph{lcm} for least common multiple, and \emph{gcd}  for greatest common divisor.
The  lcm and the  gcd  with respect to right-divisibility  of two elements  are defined in a
symmetric way.\\

An element $a\in M$ is called an \emph{atom} if $a\neq 1$ and $a=bc$ implies $b=1$ or $c=1$. A monoid~$M$ is  \emph{atomic} if $M$ is generated by its atoms and every element $g \in M$ can be expressed as a finite product of atoms. If $M$ is atomic, then  $1$ is the unique invertible element in $M$. If  $M$ is left cancellative and 1 is the unique invertible element in M, then the lcm  with respect to left-divisibility of two elements is unique, whenever it exists \cite{dehornoy}.

\subsection{Preliminaries on Garside monoids and groups}
\label{sec_bcgd_gars}
 Garside groups were introduced by P. Dehornoy and L. Paris in~\cite{DePa}, as an extension in some respects of the braid groups and the finite-type Artin groups.  These are groups of fractions of Garside monoids, and  there is the  following relation between them and other classes of monoids: Garside monoids $\subseteq$ quasi-Garside monoids  $\subseteq$ Gaussian monoids  $\subseteq$ lcm-monoids.
  \begin{defn} \cite[p.68]{dehornoy-book-braids}
 	A monoid $M$ is an \emph{lcm-monoid} if  $M$ satisfies:
 	\begin{enumerate}[(i)]
 		 		\item  $1$ is the unique invertible element in $M$ 
 		\item $M$ is  left cancellative 
	\item each pair of elements in $M$ has  a  lcm with respect to left  divisibility
 	\end{enumerate} 
 A monoid $M$ is an \emph{lcm-monoid with respect to right-divisibility} if   all the instances of left in the conditions above are replaced with right.
\end{defn}

 \begin{defn} 	
 A monoid $M$ is a \emph{Gaussian monoid} if  $M$ satisfies:
 \begin{enumerate}[(i)]
 		\item  $M$ is atomic 
 	\item $M$ is  cancellative 
 	\item each pair of elements in $M$ has  a  lcm and a gcd with respect to left and right divisibilities
 	\end{enumerate} 
 \end{defn}
A Gaussian monoid satisfies both left and right Ore's conditions, so it embeds in its group of fractions, called a  \emph{Gaussian  group}. Moreover,   every element $g \in G$ admits a unique decomposition as described in Lemma \ref{lem-unique-decomposition}.  This lemma is proved in \cite[Lemma 3.7]{dehornoy} for Garside groups and the proof holds also for Gaussian groups. 
\margin_comment{\textcolor{red}{added Lemma 1.6}}
\begin{lem}\cite{dehornoy}\label{lem-unique-decomposition}
	Let $M$ be  a  Gaussian monoid and $G$ be its group of fractions. Then every $g \in G$ admits a unique decomposition $g=uv^{-1}$, $u,v \in M$ such that $u \tilde{\wedge} v=1$, where $\tilde{\wedge}$ denotes the  gcd with respect to right divisibility. Symmetrically, $g \in G$ admits a unique decomposition $g=u'^{-1}v'$, $u',v' \in M$ such that $u' \wedge v'=1$, where $\wedge$ denotes the  gcd with respect to left  divisibility.
\end{lem}
\margin_comment{\textcolor{red}{added paragraph}}
Let $M$ be a Gaussian monoid with generating set $X$.  Let $X^*$ denote the free monoid generated by $X$.   Then $M$ has a \emph{complemented presentation}, that is $M=\operatorname{Mon}\langle X \mid R_{\lambda}\rangle$, where $R_{\lambda}$ is the set of defining relations $\{x\lambda(x,y)=y\lambda(y,x)\mid x,y \in X\}$ and $\lambda: X\times X \rightarrow X^*$ is a  function satisfying $\lambda(x,x)=1$ for every $x \in X$. The  function  $\lambda: X\times X \rightarrow X^*$ is uniquely defined, from the existence and uniqueness  of the lcm with respect to left-divisibility  in a Gaussian monoid  \cite{DePa}, \cite{dehornoy}[p.27].  \margin_comment{\textcolor{red}{added lemma 1.7}}
\begin{lem}\label{lem-lcm-right-div}
Let 	$\operatorname{Mon}\langle X \mid R_{\lambda}\rangle$ be a complemented presentation of  a Gaussian monoid $M$.  Assume the functions  $\lambda(x,-): X \rightarrow X^*$ are injective for every $x \in X$.  
 \begin{enumerate}[(i)]
 	\item Let $x \in X$.  Then  there is a   unique  $y\in X$, such that $x\lambda(x,y)=y\lambda(y,x)$  in $R_{\lambda}$.
 	\item Let $x,y \in X$, with $x \tilde{\vee}y= wx=w'y$. Then  there are  unique  $x',y'\in X$, such that $x'\lambda(x',y')$ is a word equal to $wx$  and $y'\lambda(y',x')$ is a word equal to $w'y$ in $M$. 
 \end{enumerate}
\end{lem}
\begin{proof}
	$(i)$ For every  $x \in X$,   there exists  $y\in X$, such that $x\lambda(x,y)=y\lambda(y,x)$  in $R_{\lambda}$, since 	$\operatorname{Mon}\langle X \mid R_{\lambda}\rangle$ is  a complemented presentation of  a Gaussian monoid. Assume by contradiction that there is  $z\in X$, $ z\neq y$, such that $x\lambda(x,y)=z\lambda(z,x)$ in $R_{\lambda}$.  As, by definition,  $x\lambda(x,z)=z\lambda(z,x)$, this  contradicts the injectivity of the function $\lambda(x,-)$.\\
		$(ii)$ From $x \tilde{\vee}y= wx=w'y$,  there exist $x',y'\in X$, such that	$wx=x'w_1x$ and $w'y=y'w'_1y$ in $M$. Assume by contradiction that either $x'\lambda(x',y')$ is not equal to $wx$  or $y'\lambda(y',x')$ is not  equal to $w'y$ in $M$. So, $x'\vee y'$ left divides $wx=x'w_1x$ and $w'y=y'w'_1y$, that is $wx=(x'\vee y')w_2x$ and $w'y=(x'\vee y')w'_2y$. From cancellation,  $w_2x=w'_2y$, which contradicts the minimality of $x \tilde{\vee}y= wx=w'y$. Assume by contradiction that there is an additional  $z'\in X$, such that $x'\lambda(x',z')=z'\lambda(z',x')$ is  equal to  $x'\lambda(x',y')$  in $M$. That is, 	$\lambda(x',y')=\lambda(x',z')$, which contradicts the injectivity of the function $\lambda(x',-)$.
\end{proof}
\margin_comment{\textcolor{red}{added remark}}
From Lemma \ref{lem-lcm-right-div}$(ii)$,   for every $x,y \in X$, we can assume  $R_{\lambda}$ contains $wx=w'y$, where $w,w' \in X^*$ and $wx$, $w'y$ are expressions of the lcm with respect to  right-divisibility of $x$ and $y$. 
\begin{defn} 	
	Let $M$ be  a  Gaussian monoid. Let $\Delta$ be an element in $M$.
	 \begin{enumerate}[(i)]
	 \item  $\Delta$ is a \emph{balanced element} if   its sets of left and right  divisors coincide.
		\item   $M$ is a   \emph{quasi-Garside monoid}   if  there  exists a balanced element $\Delta$  with set of  divisors,   $\Div{\Delta}$,   generating  $M$ \cite{quasi-garside}.  
		\item A  \emph{Garside monoid}  $M$ is a quasi-Garside  monoid   with   $\Div{\Delta}$  finite \cite{dehornoy}. 
			\item If   $M$ is a  Garside monoid,   $\Delta$ is called a  \emph{Garside element}. 
	\end{enumerate} 
\end{defn}
	The seminal examples of Garside groups are the braid groups and finite-type Artin groups. M. Picantin showed that torus link groups are Garside groups \cite{picantin}.  Another infinite family of  Garside groups is   the class of  structure groups  of non-degenerate  involutive set-theoretic solutions of the quantum Yang-Baxter equation \cite{chou_art}. The Garside groups obtained from these solutions satisfy many interesting properties, and we  refer to   \cite{chou_art,chou_godel1,chou_godel2}, \cite{gateva}, \cite{ deh_adv}, for more details.   In   \cite[p.70]{cedo-gars},   in the context of  these structure groups, F. Cedo shows a connection between the  Garsideness of the structure group  and the left brace associated to it.

\section{Main results}
\subsection{Construction of left $\mathscr{M}$-braces from lcm monoids}
\begin{defn}
	A \emph{left  $\mathscr{M}$-brace}  is a set  $\mathcal{B}$ with two operations, $\oplus$ and $\cdot$, such that $(\mathcal{B},\oplus)$ is a commutative  monoid, $(\mathcal{B},\cdot)$ is a monoid  and for every $a,b,c \in \mathcal{B}$,  the  axiom of left-distributivity holds:	$a \cdot (b\oplus c) = a \cdot b \oplus a \cdot c$.
	A \emph{right $\mathscr{M}$-brace} is defined similarly:  for every $a,b,c \in \mathcal{B}$,    the  axiom of right-distributivity holds: $ (a\oplus b) \cdot c  = a \cdot c \oplus b \cdot c$.\\
	$(\mathcal{B},\oplus)$  is called \emph{the  commutative monoid of the $\mathscr{M}$-brace } and $(\mathcal{B},\cdot)$  is called \emph{the  multiplicative  monoid of the $\mathscr{M}$-brace}. A \emph{two-sided  $\mathscr{M}$-brace} is a left and right $\mathscr{M}$-brace.
\end{defn}
The definition of a left $\mathscr{M}$-brace is inspired by the definition of a left brace,    with monoid structures  replacing  group structures.  However,  there is a main difference. In  a left $\mathscr{M}$-brace, the (usual) axiom of left-distributivity holds while in a left  brace this is not the usual axiom of left-distributivity, but a slightly modified one as in Equation \ref{eqn-left-dist-like}.

The definition of 	a two-sided  $\mathscr{M}$-brace is very reminiscent of that of a semiring, but there is a very crucial difference. In a semiring $M$, for every $x\in M$, $x\cdot 0=0$, where $0$ is the identity element with respect to $+$, so if  $1$ and $0$ are equal in $M$, then $M=\{0\}$, which is not the case in an  $\mathscr{M}$-brace.

\begin{proof}[Proof of Theorem \ref{thm1}]
 Let $M$ be an  lcm-monoid with respect to $\cdot$ and with identity element $1$. We define a new operation $\oplus: M\times M \rightarrow M$ in the following way. For every $a,b\in M$, we define
 \[a\oplus b =a \vee b\]
 where $\vee$ denotes the lcm of $a$ and $b$  with respect to   left divisibility. \\
 Since $M$  is an  lcm-monoid, $\oplus$ is well-defined, and for every $a,b \in M$, $a\oplus b$ belongs to $M$. Furthermore,  from \cite[p.272-273]{dehornoy},  for all $a,b,c \in M$,  $a\oplus b=b\oplus a$,   $a(b\oplus c)=ab\oplus ac$, and associativity  holds. Indeed, in the notation of \cite[p.272-273]{dehornoy}, an   lcm-monoid is a monoid that satisfies the conditions $(C_0)$, $(C_1)$ and $(C_2^+)$. 
 For  every  $a \in M$, $a\oplus1=\,a\cdot 1=\,1 \cdot a$, so $1$ is the identity  element with respect to $\oplus$ also. So,  $(M, \oplus)$ is a commutative monoid  and  $(M,\oplus,\cdot)$ is a left $\mathscr{M}$-brace.  All the elements in $M$ are idempotents with respect to $\oplus$, since for every $a \in M$, $a\oplus a=a$. \\
If  $M$ is  a Gaussian monoid with respect to $\cdot$ and with identity element $1$, then  
we define, for  every  $a,b \in M$,  $a+b= a \tilde{\vee} b$, where $\tilde{\vee} $ denotes the lcm of $a$ and $b$  with respect to  right divisibility, and  $(M,+,\cdot)$ is a right $\mathscr{M}$-brace. 
\end{proof}
In exactly the same way,  one shows that  an  lcm-monoid with respect to right-divisibility induces a right $\mathscr{M}$-brace.  Note that  $(M,\oplus,\cdot)$, as defined in the proof,  is  not a two-sided $\mathscr{M}$-brace. Indeed, if 
$ a\oplus b= a \vee b$, where $\vee$ denotes the lcm of $a$ and $b$  with respect to   left divisibility, then the  axiom of right-distributivity does not hold. 

The converse of Theorem  \ref{thm1} is not true. Indeed, given a left $\mathscr{M}$-brace,  $(\mathcal{B},\oplus, \cdot)$,  $a,b \in \mathcal{B}$,  the natural approach is to define  $ a\vee b$ by  $ a\vee b= a \oplus b$. But, as $a$ and $b$ do not necessarily left divide the element  $ a \oplus b$,  there is no meaning to this definition. The only thing we can prove is that if  for every $a \in \mathcal{B}$,  $a \oplus 1=a$, then $1$ is the unique invertible element.  Indeed, if we assume that $g \in \mathcal{B}$ is invertible, with inverse  $g^* \in \mathcal{B}$, then on one hand $g(1\oplus g ^*)=gg^*=1$, and on the second hand 
$g (1\oplus g ^*)=g \oplus gg^*=g \oplus1=g$, so  $g=1$.

\subsection{Construction of  partial braces from Gaussian groups}
In \cite{partial-s}, the authors define \emph{a partial semigroup} to be a set  $S$ together with an operation $\oplus$ that maps a subset $\mathcal{D}\subset S \times S$ into $S$ and satisfies the associative law $(g \oplus h)\oplus k=g \oplus(h\oplus k)$, in the sense that if either side is defined then so is the other and they are equal. We define a \emph{a partial monoid} to be a partial semigroup, with an identity $1$ such that $g \oplus 1$, $1\oplus g$  are always defined, equal and equal to $g$. We say that a partial monoid is \emph{commutative}, if $g \oplus h=h \oplus g$, whenever both are defined.

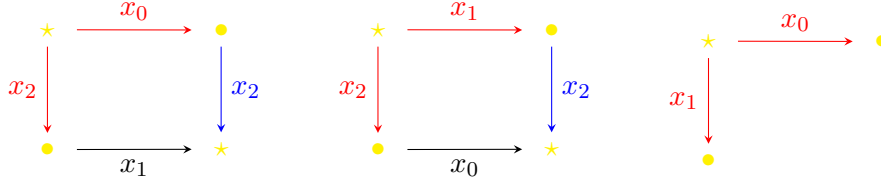
\begin{figure}[H]\label{fig-partial} 
		\begin{tikzpicture}
		\matrix (m) [matrix of math nodes,row sep=3em,column sep=4em,minimum width=2em,color=yellow]
		{ \star  & \bullet \\
			\bullet &  \star \\
		};
		\path[-stealth]
		(m-1-1) edge [red] node [above] {$x_0$} (m-1-2)
		(m-1-1) edge [red]node [left] {$x_2$} (m-2-1)
		(m-2-1) edge [] node [below] {$x_1$} (m-2-2)
		(m-1-2) edge [blue]node [right] {$x_2$} (m-2-2)
		;
		\end{tikzpicture}
		\hspace{15pt}
		\begin{tikzpicture}
		\matrix (m) [matrix of math nodes,row sep=3em,column sep=4em,minimum width=2em,color=yellow]
		{ \star  & \bullet \\
			\bullet &  \star \\
		};
		\path[-stealth]
		(m-1-1) edge [red] node [above] {$x_1$} (m-1-2)
		(m-1-1) edge [red]node [left] {$x_2$} (m-2-1)
		(m-2-1) edge [] node [below] {$x_0$} (m-2-2)
		(m-1-2) edge [blue]node [right] {$x_2$} (m-2-2)
		;
		\end{tikzpicture}
		\hspace{15pt}
		\begin{tikzpicture}
		\matrix (m) [matrix of math nodes,row sep=3em,column sep=4em,minimum width=2em,color=yellow]
		{ \star  & \bullet \\
			\bullet &  \\
		};
		\path[-stealth]
		(m-1-1) edge [red] node [above] {$x_0$} (m-1-2)
		(m-1-1) edge [red]node [left] {$x_1$} (m-2-1)
		;
		\end{tikzpicture}
		\caption{Right reversing in $M$: the mostright cube could not be completed.}
		
	\end{figure}

\begin{defn}\cite{lawson}
A \emph{regular  semigroup} is a semigroup $S$  such that for every element $s\in S$ there exists at least one   element $s^{*} \in S$ such that $ss^{*}s=s$ and $s^{*} ss^{*} =s^{*} $. The element $s^{*}$ is  called \emph{ an  inverse of $s$}.
An \emph{inverse semigroup} is a semigroup $S$  such that for every element $s\in S$ there exists a unique element $s^{*} \in S$ such that $ss^{*}s=s$ and $s^{*} ss^{*} =s^{*} $, that is $S$ is a regular semigroup such that every element in $S$ has a unique inverse. An \emph{inverse monoid}  $S$ is an inverse semigroup with multiplicative identity  $1$.
\end{defn} 
Another equivalent definition of an inverse semigroup is a regular semigroup in which all the idempotents commute.
\begin{defn}
	A \emph{partial left   brace}  is a set  $\mathcal{B}$ with two operations, $\oplus$ and $\cdot$, such that $(\mathcal{B},\oplus)$ is a commutative partial monoid, $(\mathcal{B},\cdot)$ is an inverse monoid  and for every $a,b,c \in \mathcal{B}$, such that  $b\oplus c$ and  $a \cdot b \oplus a \cdot c $ are defined, the following holds:
	\begin{equation}\label{defn-dist}
	a \cdot (b\oplus c) = a \cdot b \oplus a \cdot c 
	\end{equation}
	$(\mathcal{B},\oplus)$  is called \emph{the partial monoid of the brace } and $(\mathcal{B},\cdot)$  is called \emph{the  multiplicative inverse monoid of the brace}. 	A \emph{partial  right brace} is defined similarly:  for every $a,b,c \in \mathcal{B}$,   $ (a\oplus b) \cdot c  = a \cdot c \oplus b \cdot c$.
\end{defn}
In order to prove Theorem  \ref{thm2}, we consider a Gaussian group $G$ with respect to $\cdot$.  We define a new partial operation $\oplus$ in $G$, and in order to make the definition of $\oplus$  easier,  we use the method
of right reversing. P. Dehornoy introduced this tool in the context of braids and Garside groups. We give a very brief description of the reversing process and we refer the reader to \cite{dehornoy},  and in particular to \cite{dehornoy2} for a wider understanding of this topic.  Roughly, reversing can be used as  a  tool for constructing van Kampen diagrams in the context of presented semigroups or monoids.  As is well-known, two strings  $w$ and $w'$  are equal in  $M=\operatorname{Mon}\langle X\mid R\rangle$  if and only if there exists an $R$-derivation from $w$ to $w'$, defined to be a finite sequence of words $(w_0,..., w_p)$ such that $w_0$ is $w$, $w_p$ is $w'$, and, for each $i$, $w_{i+1}$ is obtained from $w_i$ by substituting some subword that occurs in a relation of $R$  with the other element of that relation. A van Kampen diagram for a pair of words  $(w,w')$ is a planar oriented graph with a unique source vertex and a unique  sink  vertex and edges labeled by letters of  $X$ so that the labels of each face correspond to a relation of $R$ and the labels of the bounding paths form the words $w$ and $w'$,  respectively.
The strings $w$ and $w'$ are equal in $M$ if and only if there exists a  van Kampen diagram for $(w,w')$. It  is  convenient to standardize van Kampen  diagrams,  so that they only contain vertical and horizontal edges, plus dotted arcs connecting vertices that are to be identified. Such standardized diagrams are called \emph{reversing diagrams}.  
\begin{ex}\label{ex-reversing}
	We consider the   monoid 
	$M=\operatorname{Mon} \langle X\mid  x_{1}x_{2}=x^{2}_{3}; x_{1}x_{3}=x_{2}x_{4};
	x_{2}x_{1}=x^{2}_{4}; x_{2}x_{3}=x_{3}x_{1};
	x_{1}x_{4}=x_{4}x_{2};x_{3}x_{2}=x_{4}x_{1} \rangle$. This is  a Garside monoid as it is the structure monoid of  a non-degenerate involutive set-theoretic solution of the quantum Yang-Baxter equation \cite[Ex. 1.2]{chou_aut}. 
	We illustrate with the following figure how to construct a reversing diagram, that represents a van Kampen diagram for a  pair of words $(w,w')$, such that $x_1x_2$ is the prefix of $w$ and $x_2x_1$ is the prefix of $w'$. We begin with the left-most figure, with source the star, and two paths labelled  $x_1x_2$  and $x_2x_1$ respectively.
	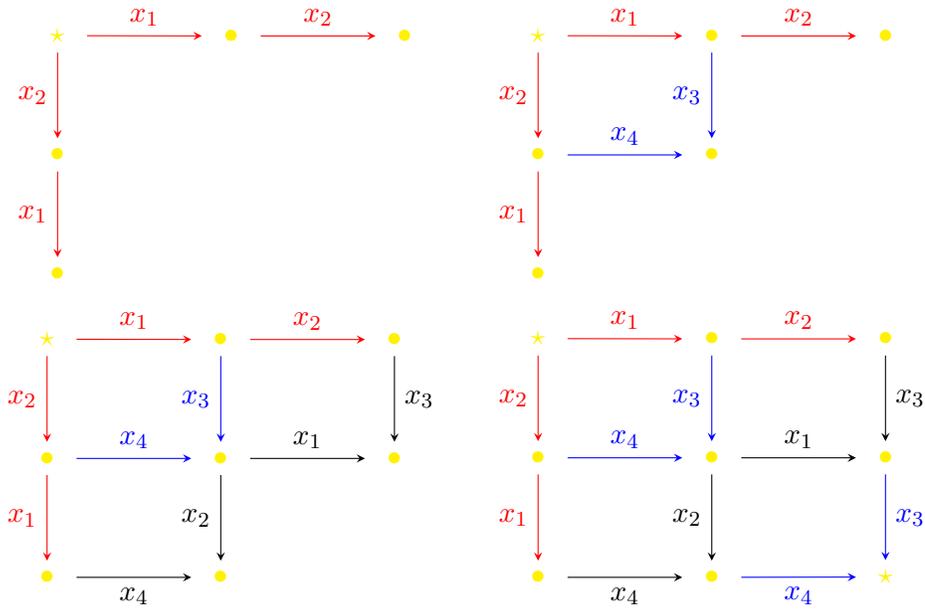
\begin{figure}[H] \label{fig-revers}
	\begin{tikzpicture}
		\matrix (m) [matrix of math nodes,row sep=3em,column sep=4em,minimum width=2em,color=yellow]
		{ 	\star & \bullet & \bullet \\
			\bullet& &  \\
			\bullet& & \\};
		\path[-stealth]
		(m-1-1) edge [red] node [above] {$x_1$} (m-1-2)
		(m-1-2) edge [red] node [above] {$x_2$} (m-1-3)
		
		(m-1-1) edge [red]node [left] {$x_2$} (m-2-1)
		(m-2-1) edge [red]node [left] {$x_1$} (m-3-1)
		
		;
		\end{tikzpicture}  
		\hspace{10pt}
	\begin{tikzpicture}
		\matrix (m) [matrix of math nodes,row sep=3em,column sep=4em,minimum width=2em,color=yellow]
		{ 	\star& \bullet & \bullet \\
			\bullet& \bullet &  \\
			\bullet& & \\};
		\path[-stealth]
		(m-1-1) edge [red] node [above] {$x_1$} (m-1-2)
		(m-1-2) edge [red] node [above] {$x_2$} (m-1-3)

		(m-1-1) edge [red]node [left] {$x_2$} (m-2-1)
		(m-2-1) edge [red]node [left] {$x_1$} (m-3-1)
		
		(m-1-2) edge [blue] node [left] {$x_3$} (m-2-2)
		(m-2-1) edge [blue] node [above] {$x_4$} (m-2-2)
		
		;
		\end{tikzpicture}

	\begin{tikzpicture}
		\matrix (m) [matrix of math nodes,row sep=3em,column sep=4em,minimum width=2em,color=yellow]
		{ 	\star& \bullet & \bullet \\
			\bullet& \bullet & \bullet \\
			\bullet& \bullet & \\};
		\path[-stealth]
		(m-1-1) edge [red] node [above] {$x_1$} (m-1-2)
		(m-1-2) edge [red] node [above] {$x_2$} (m-1-3)

		(m-1-1) edge [red]node [left] {$x_2$} (m-2-1)
		(m-2-1) edge [red]node [left] {$x_1$} (m-3-1)
		
		(m-1-2) edge [blue] node [left] {$x_3$} (m-2-2)
		(m-2-1) edge [blue] node [above] {$x_4$} (m-2-2)
		
		(m-2-2) edge [] node [above] {$x_1$} (m-2-3)
		(m-1-3) edge [] node [right] {$x_3$} (m-2-3)
		(m-2-2) edge [] node [left] {$x_2$} (m-3-2)
		(m-3-1) edge [] node [below] {$x_4$} (m-3-2)

		;
		\end{tikzpicture} 
		\hspace{10pt}
	\begin{tikzpicture}
		\matrix (m) [matrix of math nodes,row sep=3em,column sep=4em,minimum width=2em,color=yellow]
		{ 	\star& \bullet & \bullet \\
			\bullet& \bullet & \bullet \\
			\bullet& \bullet & \star\\};
		\path[-stealth]
		(m-1-1) edge [red] node [above] {$x_1$} (m-1-2)
		(m-1-2) edge [red] node [above] {$x_2$} (m-1-3)

		(m-1-1) edge [red]node [left] {$x_2$} (m-2-1)
		(m-2-1) edge [red]node [left] {$x_1$} (m-3-1)
		
		(m-1-2) edge [blue] node [left] {$x_3$} (m-2-2)
		(m-2-1) edge [blue] node [above] {$x_4$} (m-2-2)
		
		(m-2-2) edge [] node [above] {$x_1$} (m-2-3)
		(m-1-3) edge [] node [right] {$x_3$} (m-2-3)
		(m-2-2) edge [] node [left] {$x_2$} (m-3-2)
		(m-3-1) edge [] node [below] {$x_4$} (m-3-2)
		
		(m-2-3) edge [blue] node [right] {$x_3$} (m-3-3)
		(m-3-2) edge [blue] node [below] {$x_4$} (m-3-3)
		
		;
		\end{tikzpicture} 
		\caption{Right reversing  to build a van Kampen diagram for the pair $(x_1x_2x_3^2, x_2x_1x_4^2)$. }
	\end{figure} 
\end{ex}
The process of reversing is not successful for  every monoid, and for every pair of elements.  However,  if a  monoid is a Garside monoid, as in the case of the monoid $M$ in Example \ref{ex-reversing}, then right and left reversing are successful for every pair of elements in the monoid  \cite{dehornoy,dehornoy2}. But in an arbitrary monoid, it is not necessarily true anymore, and it may occur that the process never terminates, and even if it terminates    there may be  some  obstructions. Assume we have a subdiagram with horizontal edge labelled $x$ and vertical edge labelled $y$.  If there is no relation $x...=y... $ in $R$, then the subdiagram  cannot be completed, and so the diagram neither.   On the opposite, if there are more than one relation $x...=y... $ in $R$, then there may be several  different  ways to close the diagram. 
We illustrate  an obstruction of the first kind  in the use of right reversing that occurs  
	for  $M=\operatorname{Mon}\langle x_0,x_1,x_2 \mid 	 x_0x_2 =x_2x_1; x_1x_2=x_2x_0\rangle$.  As there is no defining relation $x_0...=x_1...$, the mostright  diagram in Figure \ref{fig-partial} cannot be completed.

 \begin{proof}[Proof of Theorem \ref{thm2}]
 	Let $M$ be a Gaussian monoid with respect to $\cdot$,  with generating set $X$ and identity element $1$.  Let $X^*$ denote the free monoid generated by $X$.  Let $M=\operatorname{Mon}\langle X \mid R_{\lambda}\rangle$ be  its  complemented presentation.   Let  $G$ be the  group of fractions of $M$,  $G=MM^{-1}$, where $M^{-1}=\{w^{-1}\mid w \in M\}$.  We define in $G$   a partial operation,  $\oplus$, and show that  $(G,\oplus,\cdot)$ is a partial brace.  Clearly, as 
 	$(G, \cdot)$ is a group, it is an inverse monoid. It remains to show that  $(G, \oplus)$ is a partial commutative monoid. The proof contains two parts: the iterative definition of $\oplus$ in $G$ and the proof $\oplus$ is well-defined and satisfies the relevant properties. 
 	For every $g \in G$, we define 
 	\begin{gather}\label{eq-plus}
 	g \oplus 1=g1=1g=g\\
 	g\oplus g=g 1=1g \label{eq-plusss}
 	\end{gather}
 	\margin_comment{\textcolor{red}{defn on generators}}
 		First, we define $\oplus$  for pairs of elements both in  $X$.  For every $x,y\in X$, we define:
 	\[x\oplus y=x \vee y\]
 	where $\vee$ denotes the lcm of $x$ and $y$  with respect to   left divisibility.  Clearly, $x\oplus y= x\lambda(x,y)$ and $y  \oplus x=y\lambda(y,x)$, and as these are equal in $M$, $x \oplus y=y\oplus x$.  Moreover, for each  pair $(x,y) \in X^2$, there is a unique expression  $x\lambda(x,y)$, since  the function $\lambda(x,-):X \rightarrow X^*$  is injective (see Lemma \ref{lem-lcm-right-div}$(i)$). 
 		\margin_comment{\textcolor{red}{extended defn. to pos. words }}The definition of $\oplus$  on  pairs of elements from  $X$ is extended  to pairs of elements from   $X^*$  using the reversing diagrams. Indeed, for every $a,b \in X^*$, $a \oplus b$ is always  (uniquely) defined and  $a \oplus b=a \vee b=b\oplus a$. The operation $\oplus$ is well-defined on $X^* \times X^*$. Indeed, if  $a=a'$ in $M$, where $a' \in X^*$,  that is $a$ and $a'$  are different words in $X^*$ representing the same element in $M$, then from the uniqueness of the lcm w.r to left-divisibility $a \oplus b=a' \oplus b$ in $ M$. \\
 		\margin_comment{\textcolor{red}{defn on inverses of gen.}}
 	Next, we define $\oplus$  for pairs of elements in $X^{-1} \times X^{-1}$.  From Lemma \ref{lem-lcm-right-div},   for every $x,y \in X$, we can assume  $R_{\lambda}$ contains $wx=w'y$, where $w,w' \in X^*$ and $wx$, $w'y$ are words in $X^*$ representing  the lcm with respect to  right-divisibility of $x$ and $y$. 	So, we define
 		\begin{gather}\label{eq-minus1}
 x^{-1}\oplus y^{-1}= x^{-1}w^{-1}\\
 y^{-1}\oplus x^{-1}= y^{-1}w'^{-1}
 \label{eq-minus2}	\end{gather}
 As, $x^{-1}w^{-1}=y^{-1}w'^{-1}$ in $G$,  $x^{-1}\oplus y^{-1}=  y^{-1}\oplus x^{-1}$.  	\margin_comment{\textcolor{red}{extended defn. to neg. words }}The definition of $\oplus$  on  pairs of elements from  $X^{-1}$ is extended  to pairs of elements from   $(X^{-1})^*$ iteratively, using an extended version of a  reversing diagram  in which we allow negative labelled arrows. For every $u,v \in X^*$,  $u^{-1} \oplus v^{-1}$ is always defined, and it holds that if   $a \oplus b= au=bv$, where $a,b \in X^*$,  then  $u^{-1} \oplus v^{-1}= u^{-1}  a^{-1}= v^{-1} b^{-1}=v^{-1} \oplus u^{-1}$.  The operation $\oplus$ is well-defined on $(X^{-1})^* \times (X^{-1})^*$. Indeed, if    $u=u'$ in $M$, where $u' \in X^*$, then from the uniqueness of the lcm w.r to right-divisibility,  $u' \tilde{\vee} v= u\tilde{\vee} v$, so $u^{-1} \oplus v^{-1}=u'^{-1} \oplus v^{-1}$. 
 Note that cancellativity and the  existence and uniqueness  of the lcm with respect to  left and right divisibilities  in a Gaussian monoid imply that   $\oplus$ is always  uniquely defined for pairs of elements either both in $M$ or  both in  $M^{-1}$.\\
 	 	\margin_comment{\textcolor{red}{tricky part of defn}}
  Now, we 	arrive to the definition of  $\oplus$ for pairs of elements in $X\times X^{-1}$ and $X^{-1}\times X$,  which requires far  more  carefulness.  Let $x,y\in X$,  with $x\oplus y =x u=yv$, where $u,v \in  X^*$.  If both $u$ and $v$ belong to $X$, then we define $	x^{-1} \oplus u$ and  $u \oplus x^{-1}$,  as  in  (\ref{eq-+-minus1}), 
$y^{-1} \oplus v$ and  $v	\oplus 	y^{-1}$ as in  (\ref{eq-+-minus2}), otherwise they are not defined:	\margin_comment{\textcolor{red}{restrictive defn}}
  		\begin{gather}\label{eq-+-minus1}
  	x^{-1} \oplus u=\,	x^{-1} y  \;\;\;\;\;\;\;\;\;\;\;\;    u \oplus x^{-1}=\,	uv^{-1} \\  	
  		y^{-1} \oplus v=\,	y^{-1} x \;\;\;\;\;\;\;\;\;\;\;\; v	\oplus 	y^{-1} =\,	vu^{-1}\label{eq-+-minus2}
  		\end{gather}
  	
  Clearly, 	$x^{-1} \oplus u=u \oplus x^{-1}$ and  $y^{-1} \oplus v=v	\oplus 	y^{-1}$ in $G$, whenever they are defined. If,  for $x,u,y,v\in X$,  	$x^{-1} \oplus u$  (or $y^{-1} \oplus v$) is not defined, we call $xu$ (or $yv$) a \emph{forbidden word}.  If,  for $x\in X$, there exists  $z\in X$, $z\neq x$,  such that $x\vee z=xx$, then $x \oplus x^{-1} $ is defined and $x \oplus x^{-1}\neq 1$, otherwise   $xx$ is a forbidden word and $x \oplus x^{-1} $ is not defined. In Figure \ref{fig-defn-+-}, we illustrate diagrammatically  the equations (\ref{eq-minus1})-(\ref{eq-+-minus2}) and the extended use of a reversing diagram. 
 \begin{figure}[H]
 \scalebox{0.8}[0.9]{
\begin{tikzpicture}
 	\matrix (m) [matrix of math nodes,row sep=3em,column sep=4em,minimum width=2em,color=blue,ampersand replacement=\&]
 	{ \star \& \bullet\\
 		\bullet \& \star\\};
 	\path[-stealth]
 	(m-1-1) edge [double,red] node [above] {$y$} (m-1-2)
 	
 	(m-1-1) edge [double,red]node [left] {$x$} (m-2-1)
 	
 	(m-1-2) edge[blue] node [right] {$v$} (m-2-2)
 	
 	(m-2-1) edge node [below] {$u$} (m-2-2)
 	;
 	\end{tikzpicture}}
 		\hspace{5pt}
 	 \scalebox{0.8}[0.9]{	\begin{tikzpicture}
 	\matrix (m) [matrix of math nodes,row sep=3em,column sep=4em,minimum width=2em,color=blue,ampersand replacement=\&]
 	{ \bullet \& \star\\
 	 \star 	\& \bullet\\};
 	\path[-stealth]
 	(m-1-1) edge [double,red] node [above] {$y$} (m-1-2)
 	
 	(m-2-1) edge [dotted,red]node [left] {$x^{-1}$} (m-1-1)
 	
 	(m-2-2) edge[dotted,blue] node [right] {$v^{-1}$} (m-1-2)
 	
 	(m-2-1) edge node [below] {$u$} (m-2-2)
 	;
 	\end{tikzpicture}}
 	 	\hspace{5pt}
 	  \scalebox{0.8}[0.9]{	\begin{tikzpicture}
 	\matrix (m) [matrix of math nodes,row sep=3em,column sep=4em,minimum width=2em,color=blue,ampersand replacement=\&]
 	{  \bullet \& \star\\
 		 \star \& \bullet\\};
 	\path[-stealth]
 	(m-1-2) edge [dotted,red] node [above] {$y^{-1}$} (m-1-1)
 	
 	(m-1-1) edge [double,red]node [left] {$x$} (m-2-1)
 	
 	(m-1-2) edge[blue] node [right] {$v$} (m-2-2)
 	
 	(m-2-2) edge[dotted] node [below] {$u^{-1}$} (m-2-1)
 	;
 	\end{tikzpicture}}
 		\hspace{5pt}
  \scalebox{0.8}[0.9]{	\begin{tikzpicture}
 	\matrix (m) [matrix of math nodes,row sep=3em,column sep=4em,minimum width=2em,color=blue,ampersand replacement=\&]
 	{ \star \&  \bullet  \\
 	\bullet \&	\star  \\};
 	\path[-stealth]
 	(m-1-2) edge [dotted,red] node [above] {$y^{-1}$} (m-1-1)
 	
 	(m-2-1) edge [dotted,red]node [left] {$x{-1}$} (m-1-1)
 	
 	(m-2-2) edge[dotted,blue] node [right] {$v^{-1}$} (m-1-2)
 	
 	(m-2-2) edge[dotted] node [below] {$u^{-1}$} (m-2-1)
 	;
 	\end{tikzpicture}}
\caption{From $xu=yv$ in  the reversing diagram for $x\oplus y$ (most-left), we define $x^{-1} \oplus u$,  $y^{-1} \oplus v$, $u^{-1} \oplus v^{-1}$ by reversing parallel arrows. The directed paths are read from  a star to another.}	\label{fig-defn-+-}
 \end{figure}
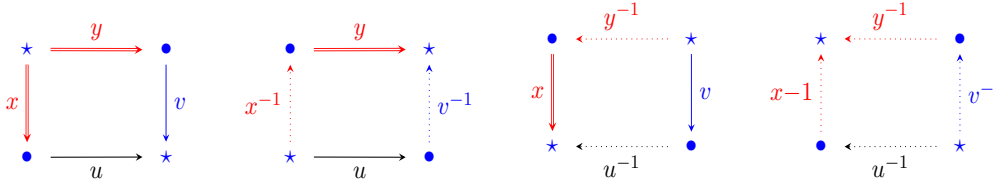
  The definition of $\oplus$  on elements
from    $X \times X^{-1}$ and $X^{-1} \times X$ is extended to  elements  in  $X^* \times (X^{-1})^*$ and $ (X^{-1})^*\times X^* $. We  illustrate  in Figure  \ref{fig-defn-oplus-mixed} the iterative computation of $u\oplus v^{-1}$, with   $u\in X^*$,  $v^{-1} \in  (X^{-1})^*$, $u=x_{i_1} x_{i_2}...x_{i_k}$, $v^{-1}=x_{j_1}^{-1}...x_{j_l}^{-1}$. From the very restrictive definition of   $\oplus$   in  $X \times X^{-1}$, and the assumption on $\lambda$,  at the completion  of the  squares (whenever possible) each  edge is labelled by a unique element from  $X \cup X^{-1}$.  So, $u\oplus v^{-1}=uz^{-1}=v^{-1 }w$, where     $w=x_{r_1}x_{r_2}...x_{r_k}$, $z^{-1}=x_{s_1}^{-1}...x_{s_l}^{-1}$. \\
	\margin_comment{\textcolor{red}{long proof  well-defined}}
We show that $\oplus$ is well-defined for elements in $M \times M^{-1}$ or  $M^{-1}\times M$.  Assume $u=u'$ in $M$, where $u' \in X^*$. Note that $v^{-1}\oplus u$ is defined if and only  if in any  word equal to $vu$  in $M$  there is no forbidden subword.    As  $vu=vu'$ in $M$, $v^{-1}\oplus u$ is defined if and only  if $v^{-1}\oplus u'$ is defined. We show that $u\oplus v^{-1}=uz^{-1}=v^{-1 }w$ and $u'\oplus v^{-1}=u'z'^{-1}=v^{-1 }w'$ are equal in $G$, whenever both are defined. From the completion of the first row in the diagram (as in Figure \ref{fig-defn-oplus-mixed}),  $x_{j_1}u=w_1x_{s_1}$  and $x_{j_1}u'=w'_1x_{s'_1}$ in $M$,   with $w_1=x_{m_1}x_{m_2}...x_{m_k}$, $w'_1=x_{m'_1}x_{m'_2}...x_{m'_k}$. From  $u=u'$ in $M$, we have  $w_1x_{s_1}=w'_1x_{s'_1}$ in $M$. So,  
 $x_{s_1}\tilde{\vee}x_{s'_1}=tx_{s_1}$, $t \in M$,  right divides $w_1x_{s_1}$ and $w'_1x_{s'_1}$, that is $w_1x_{s_1}=w_2tx_{s_1}$ and $w'_1x_{s'_1}=w'_2tx_{s_1}$. From   cancellativity in $M$, $w_2=w'_2$ and so $w_1=w'_1$, which implies $x_{s_1}=x_{s'_1}$.  Using the same argument in the completion of the diagram row by row in the same way, we obtain $w=w'$ in $M$ and $z=z'$ in $X^*$. 
 	\margin_comment{\textcolor{red}{added important diag.}}
   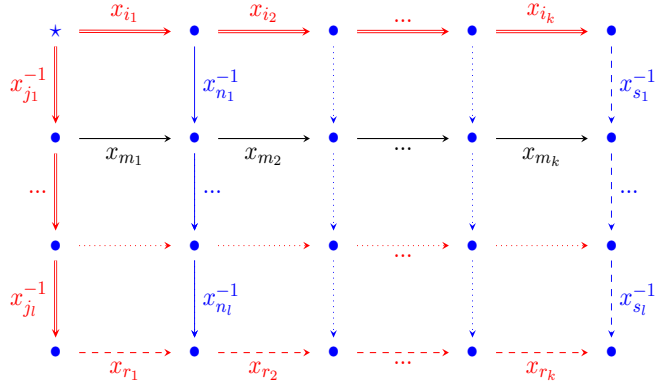
\begin{figure}[H]
  	\scalebox{0.8}[0.9]{
  		\begin{tikzpicture}
  		\matrix (m) [matrix of math nodes,row sep=3em,column sep=4em,minimum width=2em,color=blue,ampersand replacement=\&]
  		{ \star \& \bullet \& \bullet \& \bullet \& \bullet\\
  			\bullet \&  \bullet\& \bullet\& \bullet\& \bullet\\
  				\bullet \&  \bullet\& \bullet\& \bullet\& \bullet\\
  			\bullet \&  \bullet\& \bullet\& \bullet\& \bullet\\};
  		\path[-stealth]
  		(m-1-1) edge [double,red] node [above] {$x_{i_1}$} (m-1-2)
  			(m-1-2) edge [double,red] node [above] {$x_{i_2}$} (m-1-3)
  				(m-1-3) edge [double,red] node [above] {$...$} (m-1-4)
  					(m-1-4) edge [double,red] node [above] {$x_{i_k}$} (m-1-5)
  					
  		(m-1-1) edge [double,red]node [left] {$x_{j_1}^{-1}$} (m-2-1)
  		(m-2-1) edge [double,red]node [left] {$...$} (m-3-1)
  			(m-3-1) edge [double,red]node [left] {$x_{j_l}^{-1}$} (m-4-1)
  		
  		(m-1-2) edge[blue] node [right] {$x_{n_1}^{-1}$} (m-2-2)
  			(m-2-2) edge[blue] node [right] {$...$} (m-3-2)
  				(m-3-2) edge[blue] node [right] {$x_{n_l}^{-1}$} (m-4-2)
  		
  			(m-1-3) edge[dotted,blue] node [right] {} (m-2-3)
  			(m-1-4) edge[dotted,blue] node [right] {} (m-2-4)
  			(m-1-5) edge[dashed,blue] node [right] {$x_{s_1}^{-1}$} (m-2-5)

  			(m-2-2) edge[dotted,blue] node [right] {} (m-3-2)
  		(m-2-3) edge[dotted,blue] node [right] {} (m-3-3)
  		(m-2-4) edge[dotted,blue] node [right] {} (m-3-4)
  		(m-2-5) edge[dashed,blue] node [right] {$...$} (m-3-5)
  		
  		(m-3-2) edge[dotted,blue] node [right] {} (m-4-2)
  		(m-3-3) edge[dotted,blue] node [right] {} (m-4-3)
  		(m-3-4) edge[dotted,blue] node [right] {} (m-4-4)
  		(m-3-5) edge[dashed,blue] node [right] {$x_{s_l}^{-1}$} (m-4-5)
  		
  		(m-2-1) edge node [below] {$x_{m_1}$} (m-2-2)
  		(m-2-2) edge [] node [below] {$x_{m_2}$} (m-2-3)
  		(m-2-3) edge [] node [below] {$...$} (m-2-4)
  		(m-2-4) edge [] node [below] {$x_{m_k}$} (m-2-5)
  		
  			(m-3-1) edge[dotted,red] node [below] {} (m-3-2)
  		(m-3-2) edge [dotted,red] node [below] {} (m-3-3)
  		(m-3-3) edge [dotted,red] node [below] {$...$} (m-3-4)
  		(m-3-4) edge [dotted,red] node [below] {} (m-3-5)
  		
  		(m-4-1) edge [dashed,red] node [below] {$x_{r_1}$} (m-4-2)
  			(m-4-2) edge [dashed,red] node [below] {$x_{r_2}$} (m-4-3)
  		(m-4-3) edge [dashed,red] node [below] {$...$} (m-4-4)
  		(m-4-4) edge [dashed,red] node [below] {$x_{r_k}$} (m-4-5)
  		;
  		\end{tikzpicture}}
  	\caption{Computation of $u\oplus v^{-1}$, with   $u\in X^*$,  $v^{-1} \in  (X^{-1})^*$, $u=x_{i_1} x_{i_2}...x_{i_k}$, $v^{-1}=x_{j_1}^{-1}...x_{j_l}^{-1}$. }\label{fig-defn-oplus-mixed}
  	\end{figure}
  
  Assume now $v=v'$ in $M$, where $v' \in X^*$. We show that $u\oplus v^{-1}$ and $u\oplus v'^{-1}$ are equal.  From the completion of the first column  in the diagram (as in Figure \ref{fig-defn-oplus-mixed}),  $vx_{i_1}=x_{r_1}z_1$  and $v'x_{i_1}=x_{r'_1}z'_1$  in $M$,   
  with $z_1^{-1}=x^{-1}_{n_1}x^{-1}_{n_2}...x^{-1}_{n_l}$, 
   $z'^{-1}_1=x^{-1}_{n'_1}x^{-1}_{n'_2}...x^{-1}_{n'_l}$.  From  $v=v'$ in $M$, we have $x_{r_1}z_1=x_{r'_1}z'_1$ in $M$. So,  
  $x_{r_1}\vee x_{r'_1}=x_{r_1}t$, $t \in M$,  left  divides $x_{r_1}z_1$ and $x_{r'_1}z'_1$, that is $x_{r_1}z_1= x_{r_1}tz_2$ and $x_{r'_1}z'_1=x_{r_1}tz'_2$. From   cancellativity in $M$, $z_2=z'_2$ and so $z_1=z'_1$, which implies $x_{r_1}=x_{r'_1}$.  Using the same argument in the completion of the diagram column by column in the same way, we obtain $z=z'$ in $M$ and $w=w'$ in $X^*$.  So, $\oplus$ is well-defined for elements in $M \times M^{-1}$ or $M ^{-1}\times M$. Furthermore, it holds that if $v\oplus w= vu=wz$, where $u,v,w,z \in M$, and  $u \oplus v^{-1}$ is defined, then $u \oplus v^{-1}=uz^{-1}=v^{-1}w$, as it can be seen in Figure \ref{fig-defn-oplus-mixed}, when we make some changes in  the directions of the arrows as illustrated in Figure \ref{fig-defn-+-}.

	\margin_comment{\textcolor{red}{extended defn. to elements in $G$ }}  
 At last, we  define $\oplus$  for pairs of elements   in $G$.  
 	As $G$ is the group of fractions of  a Gaussian monoid $M$, every $g\in G$ can be written  uniquely as $ac^{-1}$, with $a,c\in M$, such that  $ac^{-1}$ is reduced, that is $a$ and $c$ have  no common right divisor (see lemma \ref{lem-unique-decomposition}).  Let $g=ac^{-1}$, and $h=bd^{-1}$ in $G$, $a,b,c,d\in X^*$. Again, we use the extension of the process of  reversing as a  tool  to define iteratively $g\oplus h$, as described in Figure \ref{fig-compute-inG}.
 		\begin{figure}[h] 
 	  \scalebox{0.8}[0.8]{\begin{tikzpicture}
 	 		\matrix (m) [matrix of math nodes,row sep=3em,column sep=4em,minimum width=2em,color=yellow,ampersand replacement=\&]
 		{ 	\star \& \bullet \& \bullet \\
 			\bullet \& \&  \\
 			\bullet\& \& \\};
 		\path[-stealth]
 		(m-1-1) edge [red] node [above] {$a$} (m-1-2)
 		(m-1-2) edge [dashed,red] node [above] {$c^{-1}$} (m-1-3)
 		
 		(m-1-1) edge [red]node [left] {$b$} (m-2-1)
 		(m-2-1) edge [dashed,red]node [left] {$d^{-1}$} (m-3-1)		
 		;
 		\end{tikzpicture}}
 		\hspace{20pt}
 	 \scalebox{0.8}[0.8]{\begin{tikzpicture}
 	 			\matrix (m) [matrix of math nodes,row sep=3em,column sep=4em,minimum width=2em,color=yellow,ampersand replacement=\&]
 		{ 	\star \& \bullet \& \bullet \\
 			\bullet \& \bullet \&  \\
 			\bullet\& \& \\};
 		\path[-stealth]
 			(m-1-1) edge [red] node [above] {$a$} (m-1-2)
 		(m-1-2) edge [dashed,red] node [above] {$c^{-1}$} (m-1-3)
 		
 		(m-1-1) edge [red]node [left] {$b$} (m-2-1)
 		(m-2-1) edge [dashed,red]node [left] {$d^{-1}$} (m-3-1)

 		(m-1-2) edge [blue] node [left] {$u$} (m-2-2)
 		(m-2-1) edge [] node [above] {$v$} (m-2-2) 		
 		;
 		\end{tikzpicture}}
 	 	
 		\caption{Using  diagram   to compute $g \oplus h$, where $g=ac^{-1}$,  $h=bd^{-1}$}\label{fig-compute-inG}
 	\end{figure}
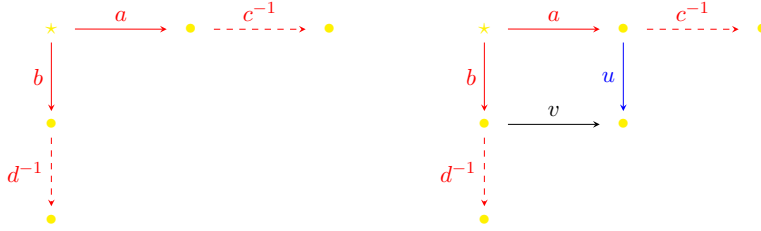 
Right reversing for the upper left (and first) square in the  diagram is always successful. But it is not necessarily the case for the subsequent squares. Indeed, from the definition of $\oplus$, $u \oplus c^{-1}$ or $v \oplus d^{-1}$ are not necessarily  defined.  If  the process of right reversing to define   $g\oplus h$ is not terminating, or if at any step, a subdiagram could not be completed, then $g \oplus h$ is not defined. If $g \oplus h$ is  defined, it can be rewritten as an element of the form $ef^{-1}$, with $e,f \in M$, since $G$ is both the group of left and right fractions of $M$.

 	We show  that 	 $\oplus$  is well-defined in $G$.   
 	Let $g,g',h \in G$,    where $g=ac^{-1}$, $g'=ef^{-1}$,  $h=bd^{-1}$  are reduced and  $a,b,c,d,e,f\in X^*$. Assume $g=g'$ in $G$. We prove that whenever $g\oplus h$ and $ g' \oplus h$ exist, then $g\oplus h$ and $g' \oplus h$ are equal in $G$.  From the uniqueness of  the decomposition of $g$ and $g'$ as fractions  in a Gaussian group, $g=g'$ in $G$ implies  $a=e$ and $c=f$ in $M$.  
 	We compute $g\oplus h$ and $ g' \oplus h$ iteratively  and show they are equal. First, at  the completion of the first square in Figure \ref{fig-well-defined-relations-inG}, we have $a\oplus b=au=bv$, $e\oplus b=eu'=bv'$, $u,u',v,v' \in X^*$. As $a=e$ in $M$,  $a\oplus b=e\oplus  b$ and from cancellation in $M$,  $u=u'$ and $v=v'$ in $M$. 
 		\margin_comment{\textcolor{red}{added  diag.}}
 		\begin{figure}[h] 
 			 \scalebox{0.8}[0.8]{\begin{tikzpicture}
 				\matrix (m) [matrix of math nodes,row sep=3em,column sep=4em,minimum width=2em,color=yellow,ampersand replacement=\&]
 				{ 	\star \& \bullet \& \bullet \\
 					\bullet \& \bullet \&  \bullet \\
 					\bullet\&  \bullet \& \star\\};
 				\path[-stealth]
 				(m-1-1) edge [red] node [above] {$a$} (m-1-2)
 				(m-1-2) edge [dashed,red] node [above] {$c^{-1}$} (m-1-3)
 				
 				(m-1-1) edge [red]node [left] {$b$} (m-2-1)
 				(m-2-1) edge [dashed,red]node [left] {$d^{-1}$} (m-3-1)

 				(m-1-2) edge [blue] node [left] {$u$} (m-2-2)
 					(m-2-1) edge [] node [above] {$v$} (m-2-2) 		
 					
 				(m-2-2) edge [dotted] node [above] {$w^{-1}_1$} (m-2-3) 		
 						(m-3-1) edge [dotted] node [below] {$w_2$} (m-3-2) 		
 								(m-3-2) edge [dotted] node [below] {$w^{-1}_3$} (m-3-3) 
 										
 			(m-2-2) edge [dotted,blue] node [left] {$z^{-1}_2$} (m-3-2)
 				(m-1-3) edge [dotted,blue] node [right] {$z_1$} (m-2-3)
 					(m-2-3) edge [dotted,blue] node [right] {$z^{-1}_3$} (m-3-3)
 									
 				;
 				\end{tikzpicture}}
 			\hspace{40pt}
 	 \scalebox{0.8}[0.8]{\begin{tikzpicture}
 		\matrix (m) [matrix of math nodes,row sep=3em,column sep=4em,minimum width=2em,color=yellow,ampersand replacement=\&]
 		{ 	\star \& \bullet \& \bullet \\
 		\bullet \& \bullet \&  \bullet \\
 		\bullet\&  \bullet \& \star\\};
 		\path[-stealth]
 		(m-1-1) edge [red] node [above] {$e$} (m-1-2)
 		(m-1-2) edge [dashed,red] node [above] {$f^{-1}$} (m-1-3)
 		
 		(m-1-1) edge [red]node [left] {$b$} (m-2-1)
 		(m-2-1) edge [dashed,red]node [left] {$d^{-1}$} (m-3-1)

 		(m-1-2) edge [blue] node [left] {$u'$} (m-2-2)
 		(m-2-1) edge [] node [above] {$v'$} (m-2-2)

 		(m-2-2) edge [dotted] node [above] {$w'^{-1}_1$} (m-2-3) 		
 		(m-3-1) edge [dotted] node [below] {$w'_2$} (m-3-2) 		
 		(m-3-2) edge [dotted] node [below] {$w'^{-1}_3$} (m-3-3) 
 		
 		(m-2-2) edge [dotted,blue] node [left] {$z^{-1}_2$} (m-3-2)
 		(m-1-3) edge [dotted,blue] node [right] {$z'_1$} (m-2-3)
 		(m-2-3) edge [dotted,blue] node [right] {$z'^{-1}_3$} (m-3-3)
 		;
 		\end{tikzpicture}}
 		\caption{Computation of $g \oplus h$, $g' \oplus h$, where $g=ac^{-1}$, $g'=ef^{-1}$,  $h=bd^{-1}$. }\label{fig-well-defined-relations-inG}
 	\end{figure}
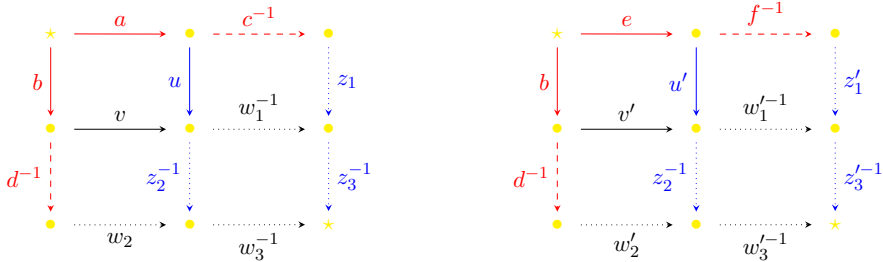 
Next,  as $c=f$, $u=u'$ in $M$, $c^{-1}\oplus u$ is defined if and only if  $f^{-1}\oplus u'$ is defined and they are equal in $G$, since $\oplus$ is well defined in $M^{-1} \times M$. The same is true for  $d^{-1}\oplus v$  and $d^{-1}\oplus v'$, since $v=v'$ in $M$. And we have $w'_1=w_1$, $z'_1=z_1$, $w'_2=w_2$ in $M$ ($z'_2=z_2$ in $X^*$). 	At last, assuming all the previous squares could be completed, it is always possible to complete the last square in Figure \ref{fig-well-defined-relations-inG}, and as $w'_1=w_1$ in $M$, we have   $z'_3=z_3$ and  $w'_3=w_3$ in $M$.   So, $g \oplus h$ is defined if and only if $g' \oplus h$ is defined and in that case  they are equal  in $G$.  Moreover, if $g \oplus h=gg'=hh'$, where $g'=z_1z_3^{-1}$ and $h'=w_2w_3^{-1}$, then  $g^{-1} \oplus g'=g^{-1}h=g'h'^{-1}$, $h^{-1} \oplus h'=h^{-1}g=h'g'^{-1}$ and  $h'^{-1} \oplus g'^{-1}=h'^{-1}h^{-1}=g'^{-1}g^{-1}$.

 	For every $g \in G$,  $a \in M$, $g$ and $gaa^{-1}$ are equal in $G$,  so, we need to  show   that whenever $g\oplus h$ and $ gaa^{-1} \oplus h$ exist, then $g\oplus h$ and $gaa^{-1} \oplus h$ are also equal in $G$. Let $g\oplus h=gg'=hh'$. Since $gaa^{-1} \oplus h$ exists,  every subdiagram can be completed, so $a\oplus g'$ exists also.  As $g' \in G$,  there exist $u,v \in M$ such that  $g'=uv^{-1}$, and we  have $a\oplus g'=ag''=g'h''$, where $g''=cd^{-1}$ and $c,d,h'' \in M$.  This implies  $a^{-1}\oplus g''=a^{-1}g'=g''h''^{-1}$, as described in Figure \ref{fig-reduction}.
 	
 		\begin{figure}[h]

 \begin{tikzpicture}[scale=0.3]
 		\matrix (m) [matrix of math nodes,row sep=3em,column sep=4em,minimum width=2em,color=yellow]
 		{ \star& \bullet&\bullet	\\
 			\bullet& \bullet&\star	\\};
 		\path[-stealth]
 		(m-1-1) edge [double,red] node [above] {$u$} (m-1-2)
 		(m-1-2) edge [double,red] node [above] {$v^{-1}$} (m-1-3)
 
 		(m-1-1) edge [double,red]node [left] {a} (m-2-1)
 		
 		(m-1-2) edge[blue] node [right] {$w$} (m-2-2)
 		(m-1-3) edge [blue]node [right] {$h''$} (m-2-3)
 	
 		(m-2-1) edge node [below] {$c$} (m-2-2)
 		(m-2-2) edge   node [below] {$d^{-1}$} (m-2-3) 	
 		;
 		\end{tikzpicture}
 		\hspace{10pt}
 	\begin{tikzpicture}[scale=0.3]
 	\matrix (m) [matrix of math nodes,row sep=3em,column sep=4em,minimum width=2em,color=yellow]
 	{ \star& \bullet&\bullet	\\
 		\bullet& \bullet&\star	\\};
 	\path[-stealth]
 	(m-1-1) edge [double,red] node [above] {$c$} (m-1-2)
 	(m-1-2) edge [double,red] node [above] {$d^{-1}$} (m-1-3)
 	
 	(m-1-1) edge [double,red]node [left] {$a^{-1}$} (m-2-1)
 	
 	(m-1-2) edge[blue] node [right] {$a^{-1}$} (m-2-2)
 	(m-1-3) edge [blue]node [right] {$h''^{-1}$} (m-2-3)
 	
 	(m-2-1) edge node [below] {$u$} (m-2-2)
 	(m-2-2) edge   node [below] {$v^{-1}$} (m-2-3)
 	;
 	\end{tikzpicture}
 		\caption{Computation of  $w\oplus g'$ at left and then  $w^{-1}\oplus g''$  at right.}
 		\label{fig-reduction}
 	\end{figure}
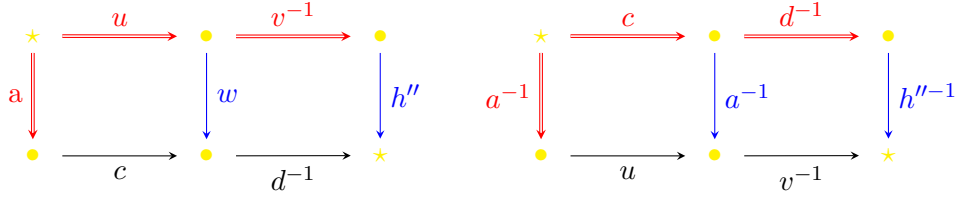
 
 	\begin{figure}[H]
 		
 			\begin{tikzpicture}
 		\matrix (m) [matrix of math nodes,row sep=3em,column sep=4em,minimum width=2em,color=yellow]
 		{ \star& \bullet\\
 			\bullet& \star\\};
 		\path[-stealth]
 		(m-1-1) edge [double,red] node [above] {$g$} (m-1-2)
 		
 		(m-1-1) edge [double,red]node [left] {$h$} (m-2-1)
 		
 		(m-1-2) edge[blue] node [right] {$g'$} (m-2-2)
 		
 		(m-2-1) edge node [below] {$h'$} (m-2-2)
 		;
 		\end{tikzpicture}
 		\hspace{20pt}
 	\begin{tikzpicture}
 		\matrix (m) [matrix of math nodes,row sep=3em,column sep=4em,minimum width=2em,color=yellow]
 		{ \star& \bullet&\bullet&\bullet	\\
 			\bullet& \bullet&\bullet&\star	\\};
 		\path[-stealth]
 		(m-1-1) edge [double,red] node [above] {$g$} (m-1-2)
 		(m-1-2) edge [double,red] node [above] {$a$} (m-1-3)
 		(m-1-3) edge [double,red] node [above] {$a^{-1}$} (m-1-4)

 		(m-1-1) edge [double,red]node [left] {h} (m-2-1)
 		
 		(m-1-2) edge[blue] node [right] {$g'$} (m-2-2)
 		(m-1-3) edge [blue]node [right] {$g''$} (m-2-3)
 		(m-1-4) edge [blue]node [right] {$g'$} (m-2-4)

 		(m-2-1) edge node [below] {$h'$} (m-2-2)
 		(m-2-2) edge   node [below] {$h''$} (m-2-3)
 		(m-2-3) edge    node [below] {$h''^{-1}$} (m-2-4)
 		;
 		\end{tikzpicture}

 		\caption{Reversing  to compute $g\oplus h$  at left and   $gaa^{-1}\oplus h$ at right.}	
 		\label{fig-additional}
 	\end{figure}
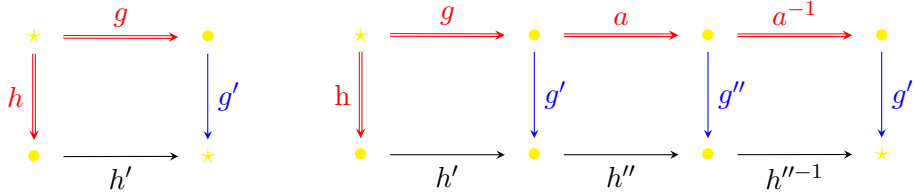
 So, we  have  $g\oplus h=gg'=hh'$ and $ gaa^{-1} \oplus h =gaa^{-1} g'=hh'h''h''^{-1}$, that is  $g\oplus h$ and $ gaa^{-1} \oplus h$ are equal in $G$.  In the same way,   $ aa^{-1}g \oplus h$ and $g \oplus h$ are  equal in $G$. So,  the partial binary operation $\oplus$ is well-defined. \\
 
  From Theorem \ref{thm1}, for every $a,b,c \in M$, $(a\oplus b) \oplus c= a \oplus(b \oplus c)$,  and  both sides are defined since these are  lcms. For $g,h \in G$, 	$g\oplus h $ (if it exists)  is not a lcm anymore, it is an element in $G$ that can be written as a  (non-trivial) word with prefix $g$ and also as a (non-trivial) word with prefix $h$, and which is computed in a very precise way, using the extension of  a reversing diagram (or a van Kampfen diagram)  as  described above. So, to compute $(g\oplus h ) \oplus k$ and $g  \oplus(h\oplus k)$, $g,h,k \in G$, is 	to find different  expressions of these elements with prefixes  $g$, $h$ and $k$. 
    	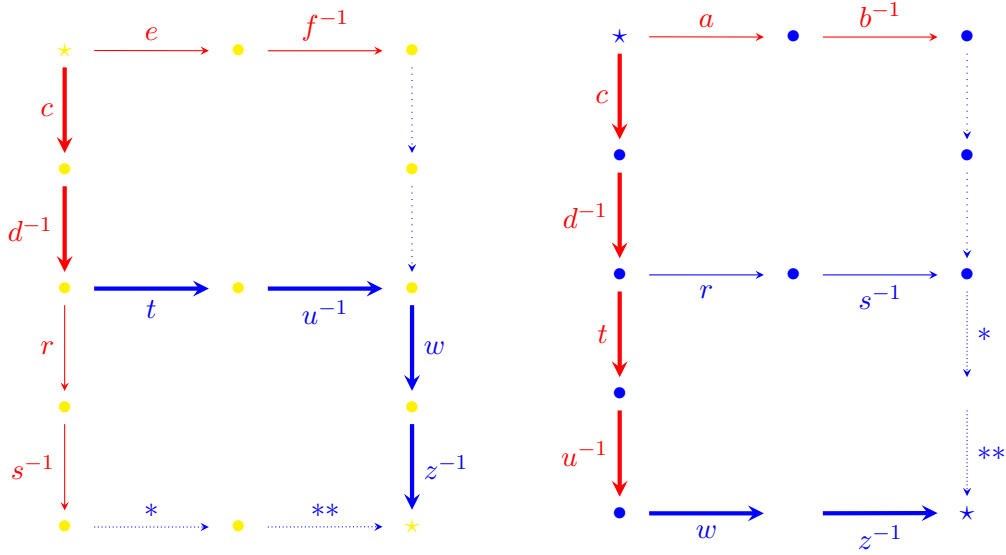
\begin{figure}[h] 	
  		\begin{minipage}{.4\textwidth}
  	\begin{flushleft}
 		\begin{tikzpicture}
 		\matrix (m) [matrix of math nodes,row sep=3em,column sep=4em,minimum width=2em,color=yellow]
 		{ 	\star & \bullet & \bullet \\
 			\bullet& &  	\bullet\\
 			\bullet & \bullet & \bullet &  \\
 			\bullet & &  \bullet\\
 			\bullet& \bullet&  \star\\};
 		\path[-stealth]
 		(m-1-1) edge [red] node [above] {$e$} (m-1-2)
 		(m-1-2) edge [red] node [above] {$f^{-1}$} (m-1-3)
 		
 		(m-1-3) edge [dotted,blue] node [right] {} (m-2-3)
 		(m-2-3) edge [dotted,blue] node [right] {} (m-3-3)
 		
 		(m-1-1) edge [ultra thick,red]node [left] {$c$} (m-2-1)
 		(m-2-1) edge [ultra thick,red]node [left] {$d^{-1}$} (m-3-1)
 		(m-3-1) edge [red]node [left] {$r$} (m-4-1)
 		(m-4-1) edge [red]node [left] {$s^{-1}$} (m-5-1)
 		
 		(m-3-1) edge [ultra thick,blue] node [below] {$t$} (m-3-2)
 		(m-3-2) edge [ultra thick,blue] node [below] {$u^{-1}$} (m-3-3)
 		
 		(m-3-3) edge [ultra thick,blue] node [right] {$w$} (m-4-3)
 		(m-4-3) edge [ultra thick,blue] node [right] {$z^{-1}$} (m-5-3)
 		
 		(m-5-1) edge [densely dotted,blue] node [above] {$*$} (m-5-2)
 		(m-5-2) edge [densely dotted,blue] node [above] {$**$} (m-5-3)
 		;
 		\end{tikzpicture}  
\end{flushleft}
 \end{minipage}
\hspace{30pt}
 \begin{minipage}{.4\textwidth}

 	\centering
 	
 		\begin{tikzpicture}
 		\matrix (m) [matrix of math nodes,row sep=3em,column sep=4em,minimum width=2em,color=blue]
 		{ 	 \star& \bullet & \bullet \\
 			\bullet& &  \bullet\\
 			\bullet&  \bullet & \bullet &  \\
 			\bullet& &  \\
 			\bullet& &  \star\\};
 		\path[-stealth]
 		(m-1-1) edge [red] node [above] {$a$} (m-1-2)
 		(m-1-2) edge [red] node [above] {$b^{-1}$} (m-1-3)
 		
 		(m-1-1) edge [ultra thick,red]node [left] {$c$} (m-2-1)
 		(m-2-1) edge [ultra thick,red]node [left] {$d^{-1}$} (m-3-1)
 		(m-3-1) edge [ultra thick,red]node [left] {$t$} (m-4-1)
 		(m-4-1) edge [ultra thick,red]node [left] {$u^{-1}$} (m-5-1)
 		
 		(m-3-1) edge [blue] node [below] {$r$} (m-3-2)
 		(m-3-2) edge [blue] node [below] {$s^{-1}$} (m-3-3)
 		
 		(m-1-3) edge [dotted,blue] node [right] {} (m-2-3)
 		(m-2-3) edge [dotted,blue] node [right] {} (m-3-3)
 		
 		(m-5-1) edge [ultra thick,blue] node [below] {$w$} (m-5-2)
 		(m-5-2) edge [ultra thick,blue] node [below] {$z^{-1}$} (m-5-3)
 		
 		(m-3-3) edge [densely dotted,blue] node [right] {$*$} (m-4-3)
 		(m-4-3) edge [densely dotted,blue] node [right] {$**$} (m-5-3)
 		
 		;
 		\end{tikzpicture}

 \end{minipage}
 		\caption{Computing $(g\oplus h ) \oplus k$ at left  an	d $g  \oplus(h\oplus k)$ at right }
\label{fig-associative}	
	\end{figure}  		
Assuming $g=ab^{-1}$, $h=cd^{-1}$, $k=ef^{-1}$  and $g\oplus h= cd^{-1}rs^{-1}$,  $h\oplus k= cd^{-1}tu^{-1}$, $a,b,c,d,e,f ,r,s,t,u\in M$,  and all the squares can be completed, we can compute $(g\oplus h ) \oplus k$ and $g  \oplus(h\oplus k)$ as described in Figure  \ref{fig-associative}. As one can see,  along the thick path from the upper left star 	to the down right star,  we read  in both diagrams a common element, so  $(g\oplus h ) \oplus k=g  \oplus(h\oplus k)$ in $G$  (whenever it is defined).  So,   $(G, \oplus)$ is a  commutative partial  monoid, with identity element equal to $1$, the identity of  $ (G,\cdot)$.

	It remains to show that Equation \ref{defn-dist} holds, that is $w(g\oplus h)=wg\oplus wh$, for every $g,h,w \in G$, such  that this expression is defined. From $g\oplus h=gg'=hh'$, multiplying at left by $w$, we have $w (g\oplus h)=wgg'=whh'$, and using the reversing diagram for $wg\oplus wh$, and Equation  \ref{eq-plusss}, we obtain $wg\oplus wh=wgg'=whh'$,  as described in Figure  \ref{fig-last}. That is,  $w (g\oplus h)=wg\oplus wh$. So, $(G, \cdot,\oplus)$ is a partial  left brace.
 	
	\begin{figure}[H] 
	\begin{tikzpicture}
	\matrix (m) [matrix of math nodes,row sep=3em,column sep=4em,minimum width=2em,color=yellow]
	{ \star& \bullet\\
		\bullet& \star\\};
	\path[-stealth]
	(m-1-1) edge [double,red] node [above] {$g$} (m-1-2)
	
	(m-1-1) edge [double,red]node [left] {$h$} (m-2-1)
	
	(m-1-2) edge[blue] node [right] {$g'$} (m-2-2)
	
	(m-2-1) edge node [below] {$h'$} (m-2-2)
	;
	\end{tikzpicture}
	\hspace{50pt}
	\begin{tikzpicture}
	\matrix (m) [matrix of math nodes,row sep=3em,column sep=4em,minimum width=2em,color=yellow]
	{ 	\star& \bullet & \bullet \\
		\bullet& \bullet & \bullet \\
		\bullet& \bullet & \star\\};
	\path[-stealth]
	(m-1-1) edge [double,red] node [above] {$w$} (m-1-2)
	(m-1-2) edge [double,red] node [above] {$g$} (m-1-3)

	(m-1-1) edge [double, red]node [left] {$w$} (m-2-1)
	(m-2-1) edge [double, red]node [left] {$h$} (m-3-1)
	
	(m-1-2) edge [blue] node [left] {$1$} (m-2-2)
	(m-2-1) edge [] node [above] {$1$} (m-2-2)
	
	(m-2-2) edge [] node [above] {$g$} (m-2-3)
	(m-1-3) edge [blue] node [right] {$1$} (m-2-3)
	(m-2-2) edge [blue] node [left] {$h$} (m-3-2)
	(m-3-1) edge [] node [below] {$1$} (m-3-2)
	
	(m-2-3) edge [blue] node [right] {$g'$} (m-3-3)
	(m-3-2) edge [] node [below] {$h'$} (m-3-3)
	
	;
	\end{tikzpicture}

	\caption{Right reversing  to compute $g \oplus h$ at left and $ wg\oplus wh$ at right}
	\label{fig-last}
\end{figure}
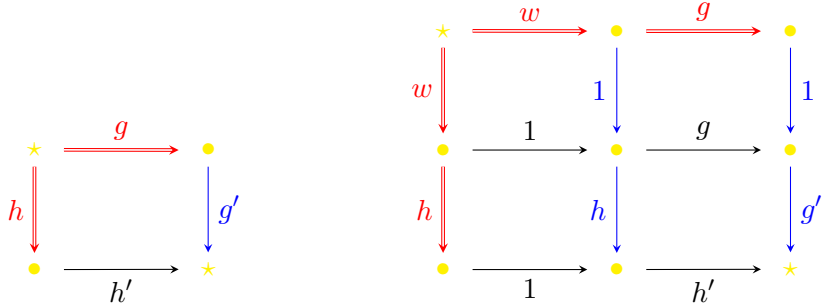 
 	\end{proof}
 	\margin_comment{\textcolor{red}{added rmk}} 
 There exists an infinite family of  Gaussian monoids satisfying  Theorem 2. Indeed, the Garside monoids arising from (non-degenerate and involutive)  set-theoretic solutions of the quantum   Yang-Baxter equation have a complemented presentation with $\lambda: X \times X \rightarrow X$ and  all the functions $\lambda(x,-): X\rightarrow X$ are injective \cite{chou_art}.  Moreover, the defining relations are length-preserving and of a simple form. Example \ref{ex-reversing} is an example of such a Garside monoid. However, not all the Garside monoids have a presentation with length-preserving defining relations, as the following example illustrates it. 
 	\margin_comment{\textcolor{red}{added example}} 
 \begin{ex}
 	Let $M=\operatorname{Mon}\langle x_1,x_2 \mid x_1x_2x_1=x_2^2\rangle$ a Garside  monoid \cite{dehornoy}.  Let $G$ be its group of fractions. From the proof of Theorem $2$, $x_1\oplus x_2=x_1(x_2x_1)=x_2(x_2)$, $x_1^{-1}\oplus x_2^{-1}=x_1^{-1}(x_2^{-1}x_1^{-1})= x_2^{-1}(x_2^{-1})$. As the definition of $\oplus$ for pairs in $X \times X^{-1}$  in the proof of Theorem $2$ is very restrictive,   it  is not possible to define  $x_i\oplus x_j^{-1}$ for $(i,j) \in \{(1,1),(1,2), (2,1),(2,2)\}$. If we try to relax this  restrictive definition, we need to add  more  rules on how to define $x_2\oplus x_2^{-1}$ in order to make $\oplus$ well defined. As an example, if we try to define $x_2\oplus x_2^{-1}$, we find there are at least two options: $x_2\oplus x_2^{-1}=x_2x_1^{-1}=x_2^{-1}(x_1x_2)$  or 
 	$x_2\oplus x_2^{-1}=x_2(x_1^{-1}x_2^{-1})=x_2^{-1}(x_1)$.  
 \end{ex}
	\margin_comment{\textcolor{red}{added details }} 
Note that even if $M$ is a Garside monoid satisfying many good properties, the operation $\oplus$ is still a partial one.  Indeed, if we consider the group of fractions  $G$ of the Garside monoid   from Example \ref{ex-reversing}, then $x_3 \oplus x_4^{-1}$, and  $x_4 \oplus x_3^{-1}$ are not defined,  since there is no $z,z'\in M$  such that $x_3\vee z=x_3x_4$,
$x_4\vee z'=x_4x_3$. Moreover, any reversing diagram requiring the definition of $x_3 \oplus x_4^{-1}$, or   $x_4 \oplus x_3^{-1}$ cannot be completed. The definition of $x_i\oplus x_j^{-1}$, $1\leq i,j \leq 4$,   is possible for all the pairs $(i,j)$, except $(1,1), (2,2), (3,4), (4,3)$. 
 \section{Appendix: Reversing diagram to show  $aa^{-1}g\oplus h$  and  $g\oplus h$ are equal in $G$}
 We use here greek letters to denote elements in $M$, and we assume $g=\alpha \beta^{-1}$, $h=\gamma \delta^{-1}$.
 \begin{figure}[H]
 	\begin{minipage}{.5\textwidth}
 		\begin{flushleft}
 			\begin{tikzpicture}[scale=0.3]
 			\matrix (m) [matrix of math nodes,row sep=3em,column sep=4em,minimum width=2em,color=yellow]
 			{ \star& \bullet&\bullet&\bullet & \bullet 	\\
 				\bullet & \bullet&\bullet & \bullet &\bullet	\\
 				\bullet & \bullet&\bullet & \bullet &\star	\\};
 			\path[-stealth]
 			(m-1-1) edge [double,red] node [above] {$a$} (m-1-2)
 			(m-1-2) edge [double,red] node [above] {$a^{-1}$} (m-1-3)
 			(m-1-3) edge [double,red] node [above] {$\alpha$} (m-1-4)
 			(m-1-4) edge [double,red] node [above] {$\beta^{-1}$} (m-1-5)

 			(m-1-1) edge [double,red]node [left] {$\gamma$} (m-2-1)
 			(m-2-1) edge [double,red]node [left] {$\delta^{-1}$} (m-3-1)
 			
 			(m-1-2) edge[blue] node [right] {$a'$} (m-2-2)
 			(m-1-3) edge [blue]node [right] {$\gamma$} (m-2-3)
 			(m-1-4) edge [blue]node [right] {$\alpha'$} (m-2-4)
 			(m-1-5) edge [blue]node [right] {$\beta'$} (m-2-5)

 			(m-2-1) edge node [below] {$\phi$} (m-2-2)
 			(m-2-2) edge   node [below] {$\phi^{-1}$} (m-2-3)
 			(m-2-3) edge    node [below] {$\gamma'$} (m-2-4)
 			(m-2-4) edge    node [below] {$\tau^{-1}$} (m-2-5)

 			(m-2-2) edge[blue] node [right] {$\psi^{-1}$} (m-3-2)
 			(m-2-3) edge [blue]node [right] {$\delta^{-1}$} (m-3-3)
 			(m-2-4) edge [blue]node [right] {$\epsilon^{-1}$} (m-3-4)
 			(m-2-5) edge [blue]node [right] {$\mu^{-1}$} (m-3-5)

 			(m-3-1) edge node [below] {$\pi$} (m-3-2)
 			(m-3-2) edge   node [below] {$\pi^{-1}$} (m-3-3)
 			(m-3-3) edge    node [below] {$\delta'$} (m-3-4)
 			(m-3-4) edge    node [below] {$\sigma^{-1}$} (m-3-5)
 			;
 			\end{tikzpicture}
 		\end{flushleft}
 	\end{minipage}
 	\hspace{75pt}
 	\begin{minipage}{.3\textwidth}
 		
 		\begin{tikzpicture}[scale=0.3]
 		\matrix (m) [matrix of math nodes,row sep=3em,column sep=4em,minimum width=2em,color=yellow]
 		{ \star& \bullet & \bullet \\
 			\bullet&  \bullet &\bullet \\
 			\bullet&  \bullet &\star\\};
 		\path[-stealth]
 		(m-1-1) edge [double,red] node [above] {$\alpha$} (m-1-2)
 		(m-1-2) edge [double,red] node [above] {$\beta^{-1}$} (m-1-3)

 		(m-1-1) edge [double,red]node [left] {$\gamma$} (m-2-1)
 		(m-2-1) edge [double,red]node [left] {$\delta^{-1}$} (m-3-1)
 		
 		(m-1-2) edge[blue] node [right] {$\alpha'$} (m-2-2)
 		(m-1-3) edge[blue] node [right] {$\beta'$} (m-2-3)
 		
 		(m-2-2) edge [blue]node [right] {$\epsilon^{-1}$} (m-3-2)
 		(m-2-3) edge[blue]node [right] {$\mu^{-1}$} (m-3-3)
 		
 		(m-2-1) edge node [below] {$\gamma'$} (m-2-2)
 		(m-3-1) edge node [below] {$\delta'$} (m-3-2)
 		
 		(m-2-2) edge node [below] {$\tau^{-1}$} (m-2-3)
 		(m-3-2) edge node [below] {$\sigma^{-1}$} (m-3-3)
 		
 		;
 		\end{tikzpicture}
 	\end{minipage}
 	\caption{Right reversing  to compute $aa^{-1}g\oplus h$ at left and  $g\oplus h$ at right.}
 	\label{fig-well-def2}
 \end{figure}
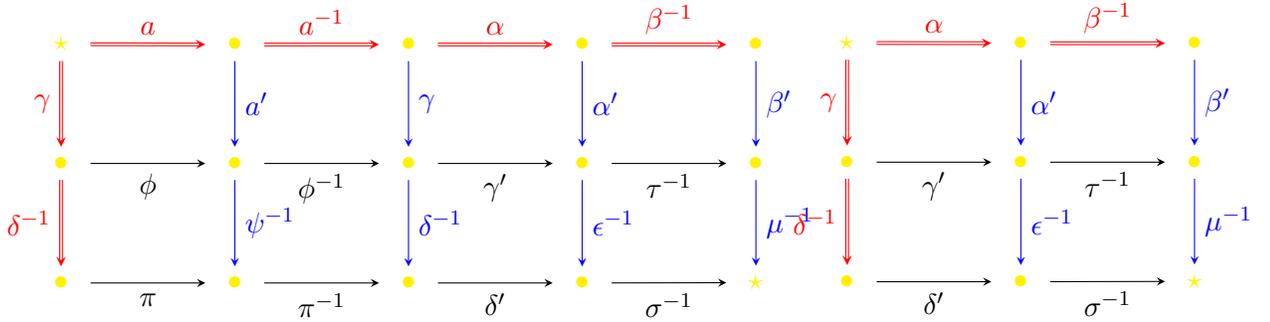

\begin{acknowledgment}
	I am very grateful to the referee for the careful reading of the paper and for his/her  comments  which helped  improve the paper.\end{acknowledgment}

\bigskip\bigskip\noindent
{ Fabienne Chouraqui}

\smallskip\noindent
University of Haifa at Oranim, Israel.

\smallskip\noindent
E-mail: {\tt fabienne.chouraqui@gmail.com} \\

{\tt fchoura@sci.haifa.ac.il}
\end{document}